\documentclass[aos,MSNbibl,rotating,seceqn,dvips]{arximspdf}
\usepackage{dcolumn}
\usepackage{graphicx}

\doi{10.1214/12-AOS1072} 
\volume{41}
\issue{1}
\pubyear{2013}
\firstpage{250}
\lastpage{268}

\makeatletter
\newcommand{\rrvert}{\vert}
\newcommand{\llvert}{\vert}
\newcolumntype{d}[1]{D{.}{.}{#1}}
\def\A{\mathbf{A}}
\def\B{\mathbf{B}}
\def\u{\mathbf{u}}
\def\x{\mathbf{x}}
\def\D{\mathbf{D}}
\def\I{\mathbf{I}}
\def\f{\mathbf{f}}
\def\g{\mathbf{g}}
\def\z{\mathbf{z}}
\def\S{\mathbf{S}}
\def\mR{\mathrm{R}}
\def\n{\nonumber}
\def\cov{\operatorname{cov}}
\def\var{\operatorname{var}}
\def\eff{_{\mathrm{eff}}}
\def\corr{\operatorname{corr}}
\def\vecl{\operatorname{vecl}}
\def\vecs{\operatorname{vec}}
\def\bSig{\bolds{\Sigma}}
\def\tY{{\widetilde{Y}}}
\def\tx{{\widetilde{\x}}}
\def\cs{\calS_{Y\mid\x}}
\def\calS{\mathcal{S}}
\def\hb{\widehat{\bb}}
\def\bb{\bolds{\beta}}
\def\ba{\bolds{\alpha}}
\def\bg{\bolds{\gamma}}
\def\defby{\stackrel{\mathrm{def}}{=}}
\newproclaim{Rem}{Remark}
\newtheorem{Th}{Theorem}

\def\cal{\mathcal}
\def\mid{|}
\makeatother

\begin{document}
\begin{frontmatter}

\title{Efficient estimation in sufficient dimension reduction}%
\runtitle{Efficient estimation in sufficient dimension reduction}

\begin{aug}
\author[A]{\fnms{Yanyuan} \snm{Ma}\thanksref{t1}\ead[label=e1]{ma@stat.tamu.edu}}
\and
\author[B]{\fnms{Liping} \snm{Zhu}\corref{}\thanksref{t2}\ead[label=e2]{zhu.liping@mail.shufe.edu.cn}}
\thankstext{t1}{Supported by NSF Grants DMS-12-06693 and DMS-10-00354 and
the National Institute of
Neurological Disorders and Stroke Grant R01-NS073671.}
\thankstext{t2}{Supported by Natural Science
Foundation of China (11071077),
Innovation Program of Shanghai Municipal Education Commission (13ZZ055),
Pujiang Project of Science and Technology Commission of Shanghai
Municipality (12PJ1403200) and Grants for New Century Excellent Talents in University,
Ministry of Education (NCET-12-0901). All the correspondence should be
directed to Liping Zhu at \printead*{e2}.}
\runauthor{Y. Ma and L. Zhu}

\affiliation{Texas A\&M University and Shanghai University of Finance
and Economics}

\address[A]{Department of Statistics\\
Texas A\&M University\\
3143 TAMU\\
College Station, Texas 77843-3143\\
USA\\
\printead{e1}} 

\address[B]{School of Statistics and Management\\
\quad and the Key Laboratory\\
\quad of Mathematical Economics\\
\quad Ministry of Eduction\\
Shanghai University of Finance and Economics\\
777 Guoding Road\\
Shanghai 200433\\
P.R. China\\
\printead{e2}}
\end{aug}

\received{\smonth{11} \syear{2011}}
\revised{\smonth{8} \syear{2012}}

%
\begin{abstract}
We develop an efficient estimation procedure for identifying and
estimating the central subspace.
Using a new way of parameterization, we convert the problem
of identifying the central subspace to the problem
of estimating a finite dimensional parameter in a semiparametric
model. This conversion allows us to derive an
efficient estimator which reaches the optimal
semiparametric efficiency bound. The resulting efficient estimator can
exhaustively estimate the central subspace
without imposing any distributional assumptions.
Our proposed efficient estimation also provides a possibility for
making inference
of parameters that uniquely identify the central subspace.
We conduct simulation studies and a real data analysis
to demonstrate the finite sample performance
in comparison with several existing
methods.\vspace*{-1pt}
\end{abstract}

%


\end{frontmatter}

\section{Introduction} \label{secintro}

Consider a general model in which
the univariate response variable $Y$ is assumed to depend on the
$p$-dimensional covariate vector $\x$ only through
a small number of linear combinations $\bb^{\mathrm{T}}\x$, where
$\bb$ is a
$p\times d$ matrix with $d<p$. In this model, how $Y$ depends on
$\bb^{\mathrm{T}}\x$ is left unspecified. It is not difficult to see that
$\bb$ is not identifiable.
The quantity of general interest is usually
the column space of $\bb$, which is
termed the central subspace if $d$ is the smallest possible value to
satisfy the model assumption~\cite{Cook1998}.

This general model was proposed by Li~\cite{Li1991} and has
attracted much
attention in the last two decades. It generated the field of
sufficient dimension reduction~\cite{Cook1998}, in which the main
interest is to estimate
the central subspace consistently.
Influential works in this area include, but are not limited to, sliced
inverse regression~\cite{Li1991},
sliced average variance estimation
\cite{CookWeisberg1991}, directional regression~\cite{LiWang2007},
the generalization of the aforementioned methods to nonelliptically
distributed predictors~\mbox{\cite{LiDong2009,DongLi2010}},
Fourier transformation~\cite{ZhuZeng2006},
cumulative slicing
estimators~\cite{ZhuZhuFeng2010} 
and conditional density based minimum average
variance estimation~\cite{Xia2007}, etc.

Despite the various estimation methods, it is unclear if any of
these estimators are optimal in the sense that they can exhaustively estimate
the entire central subspace and have the minimum
possible asymptotic estimation variance.
To the best of our knowledge, the efficiency issue has never been
discussed in the context of sufficient dimension reduction.

In this paper we study the estimation and inference in sufficient
dimension reduction. We propose a simple parameterization so that the
central subspace is uniquely identified by a $(p-d)d$-dimensional
parameter that is not subject to any constraints.
Thus we
convert the problem of identifying the central subspace into a
problem of estimating a finite dimensional parameter in a
semiparametric model.
This allows us to derive the estimation procedures and
perform inference using semiparametric tools.
How to make inference about the central subspace is a challenging issue.
This is partially caused by the complexity of estimating a space rather than
a parameter. Our new parameterization overcomes this complexity and
permits a relatively straightforward calculation of the estimation variability.

We further construct an
efficient
estimator, which reaches the minimum asymptotic estimation variance
bound among all possible consistent estimators.
Efficiency bounds are of fundamental importance to the theoretical
consideration. Such bounds quantify the minimum efficiency loss that
results from generalizing one restrictive model to a more flexible
one, and hence they can be important in making the decision of which model
to use. The efficiency bounds also provide a gold standard by which
the asymptotic efficiency of any particular semiparametric estimator
can be measured~\cite{Newey1990}. Generally speaking, a
semiparametric efficient estimator is usually the ultimate destination
when searching for consistent estimators or trying to
improve existing procedures.
When an efficient estimator is
obtained, the procedure of estimation can be considered to have
reached certain optimality.\looseness=-1

In the literature, vast and significant effort has
been devoted to studying the semiparametric efficiency bounds for
consistent estimators in semiparametric models.
The simplest and
most familiar examples are the ordinary and weighted least square
estimators in the linear regression setting. Efficiency issues are also
considered
in more complex semiparametric problems such as
regressions with missing\vadjust{\goodbreak} covariates~\cite{RobinsRotnitzkyZhao1994},
skewed distribution families~\cite{MaGentonTsiatis2005,MaHart2007},
measurement error models~\cite{TsiatisMa2004,MaCarroll2006},
partially linear models~\cite{MaChiouWang2006}, the Cox model
\cite{Tsiatis2006}, page 113, accelerated failure
model~\cite{ZengLin2007a} or other
general survival models~\cite{ZengLin2007b} and
latent variable models~\cite{MaGenton2010}.

One typical semiparametric tool is to obtain estimators through
obtaining the corresponding influence functions.
In deriving the influence function family and its efficient member, we
use the geometric technique illustrated in~\cite{BickelKlaassenRitovWellner1993} and
\cite{Tsiatis2006}. All our derivations are performed
without using the linearity or constant variance
condition that is often assumed in the dimension reduction
literature. Our analysis is thus readily applicable when some
covariates are discrete or categorical. In summary, we provide an
efficient estimator which can exhaustively estimate the
central subspace without imposing any distributional assumptions on the
covariate~$\x$.

The rest of this paper is organized as follows. In Section
\ref{secsemi}, we propose a simple parameterization of the central
subspace and
highlight the semiparametric approach to estimating the central subspace.
We also derive the efficient score function. In Section~\ref{seceff}, we present a class of locally efficient estimators and
identify the efficient member. We
illustrate how to implement the efficient estimator to reach the
optimal efficiency bound. Simulation studies are conducted in Section~\ref{secsimu} to demonstrate the finite sample performance and the
method is implemented in a real data example in Section~\ref{secexample}.
We finish the paper with a brief discussion in Section~\ref{secdiscuss}.
All the technical derivations are given in a supplementary material~\cite{supp}.

\section{The semiparametric formulation}\label{secsemi}
\subsection{Parameterization of central subspace}
In the context of sufficient dimension reduction
\cite{Li1991,Cook1998}, one often assumes
%
\begin{equation}
\label{eqcs} F(y\mid\x) = F\bigl(y\mid\bb^{\mathrm{T}}\x\bigr)\qquad \mbox{for }
y\in \mR,
\end{equation}
where
$F(y\mid\x) \defby\operatorname{Pr}(Y\le
y\mid\x)$ is the conditional distribution function of
the response $Y$ given the covariates $\x$, and $\bb$ is a $p\times d$
matrix as defined previously. The goal of sufficient dimension
reduction is to estimate the column space of $\bb$, which is termed
the dimension reduction subspace.
Because a dimension reduction subspace is not necessarily unique,
the primary interest is usually the central subspace $\cs$, which is
defined as
the minimum dimension reduction subspace if it exists and is unique
\cite{Cook1998}. The dimension of $\cs$, denoted with $d$, is
commonly referred to as the
structural dimension.
Similarly to~\cite{Cook1994},
we exclude a pathological case
where there exists a vector $\ba$ such that
$\ba^{\mathrm{T}}\x$
is a deterministic function of
$\bb^{\mathrm{T}}\x$ while $\ba$ does not belong to
the column space of $\bb$.

The central subspace $\cs$ has a well-known invariance property~\cite{Cook1998},\break page~106, that is,
$\cs= \D\calS_{Y\mid{\z}}$, where $\z= \D^{\mathrm{T}}\x+
\mathbf{b}$
for any $p\times p$ nonsingular matrix $\D$ and any length $p$ vector
$\mathbf{b}$.
This allows us to assume throughout that\vadjust{\goodbreak} the covariate vector $\x$
satisfies $E(\x)=\mathbf{0}$ and $\cov(\x)=\I_p$.
Identifying $\cs$ is the essential interest of sufficient dimension
reduction for
model (\ref{eqcs}). Typically, $\cs$ is identified through estimating
a basis matrix $\bb\in
\mR^{p\times d}$ of minimal dimension that satisfies~(\ref{eqcs}). Although
$\cs$ is unique, the basis matrix $\bb$ is clearly not. In fact, for
any $d\times d$ full rank
matrix $\mathbf{ A}$, $\bb\mathbf { A}$ generates the same column space as $\bb$.
Thus, to uniquely map one central subspace $\cs$ to one basis matrix,
we need to focus on one representative member of all the $\bb\A$
matrices generated by different $\A$'s.
We write $\bb=(\bb_u^{\mathrm{T}}, \bb_l^{\mathrm{T}})^{\mathrm{T}}$,
where the upper submatrix $\bb_u$ has size $d\times d$ and the lower
submatrix $\bb_l$ has size $(p-d)\times d$.
Because $\bb$ has rank $d$, we can assume without loss of generality
that $\bb_u$ is invertible. The advantage of using $\bb\bb_u^{-1}$ is
that its upper $d\times d$ submatrix is the identity matrix, while the
lower $(p-d)\times d$ matrix can be any matrix. In addition, two
matrices $\bb_1\bb_{1u}^{-1}$ and $\bb_2\bb_{2u}^{-1}$ are different
if and only if the column spaces of $\bb_1$ and $\bb_2$ are
different. Therefore, if we consider the set of all the $p\times d$
matrices $\bb$ where the upper $d\times d$ submatrix is the identity
matrix $\I_d$, it has a one-to-one mapping with the set of all the
different central subspaces. Thus, as long as we restrict our attention
to the set of all such matrices,
the problem of identifying $\cs$ is converted to the problem of
estimating $\bb_l$, which contains $p_t=(p-d)d$
free parameters. Note that $p_t$ is the dimension of
the Grassmann manifold formed by the column spaces of all different
$\bb$ matrices. Thus, we can view $\bb_l$ as a unique parameterization
of the manifold.
Here the subscript ``$_t$'' stands for total.
For notational convenience in the remainder of the text, for an
arbitrary $p\times d$ matrix $\bb=(\bb_u^{\mathrm{T}}, \bb_l^{\mathrm{T}})^{\mathrm{T}}$,
we define the concatenation of the columns contained in the
lower $p-d$ rows of $\bb$ as
$\vecl(\bb)=\vecs(\bb_l) = (\beta_{d+1,1}, \ldots, \beta_{p,1},
\ldots,
\beta_{d+1,d}, \ldots, \beta_{p,d})^{\mathrm{T}}$, where in the
notation $\vecl$, ``vec''
stands for vectorization, and ``l'' stands for the lower part of the
original matrix.
We then can write the
concatenation of the parameters in $\bb$ as
$\vecl(\bb)$. Thus, from now on, we only consider basis matrix of
$\cs$ that has the form $\bb=(\I_d, \bb_l^{\mathrm{T}})^{\mathrm
{T}}$, where
$\bb_l$
is a $(p-d)\times d$ matrix. Estimating the parameters in $\bb$ is a typical
semiparametric estimation problem, in which the parameter of interest is
$\vecl(\bb)$.
Therefore we have converted the problem of estimating the central
space $\cs$ into a problem of semiparametric estimation.

\begin{Rem}
The above parameterization of $\cs$ excludes the
pathological case where one or more of the first $d$ covariates
do not contribute to the model or contribute to the model through a
fixed linear combination. When this happens,
$\bb_u$ will be singular. However,
because $\bb$ has rank $d$, hence if this
happens, one can always rotate the order of
the covariates (hence rotate the rows of $\bb$) to ensure that
after rotation, the resulting $\bb_u$ has full rank.
\end{Rem}

\subsection{Efficient score}
In this section we derive the efficient score
for estimating $\bb$ under the above parameterization.
That is, we now consider model (\ref{eqcs}), where $\bb=(\I_d,
\bb_l^{\mathrm{T}})^{\mathrm{T}}$ and\vadjust{\goodbreak} $\x$ satisfies $E(\x)=\mathbf{0}$ and
$\var(\x)=\I_p$.
The general semiparametric technique we use is originated from
\cite{BickelKlaassenRitovWellner1993} and is wonderfully
presented in~\cite{Tsiatis2006}.
Using this approach, we obtain the main result of this section,
that we can use (\ref{eqcseff})
to obtain an efficient estimation of~$\bb$.\looseness=-1

The likelihood of one random observation
$(\x,Y)$ in (\ref{eqcs}) is $\eta_1(\x)\eta_2(Y,\bb^{\mathrm
{T}}\x)$,
where
$\eta_1$ is a probability mass function (p.m.f.) or a probability density
function (p.d.f.) of
$\x$, or a mixture, depending on whether $\x$ contains
discrete variables, and $\eta_2$ is the conditional p.m.f./p.d.f.
of $Y$ on $\x$. We view
$\eta_1, \eta_2$ as infinite dimensional nuisance parameters and
$\vecl(\bb)$ as the $p_t$-dimensional parameter of interest.
Following the semiparametric analysis procedure, we first derive the
nuisance tangent space $\Lambda=\Lambda_1\oplus\Lambda_2$, where
\begin{eqnarray*}
\Lambda_1&=&\bigl\{\f(\x)\dvtx\forall\f\mbox{ such that } E(\f
)=\mathbf{ 0}\bigr\},
\\
\Lambda_2&=&\bigl\{\f\bigl(Y,\bb^{\mathrm{T}}\x\bigr)\dvtx\forall\f
\mbox{ such that } E(\f\mid \x)=E\bigl(\f\mid\bb^{\mathrm{T}}\x\bigr)=\mathbf{0}
\bigr\}.
\end{eqnarray*}
Here, the notation $\oplus$ means the
usual addition of the two spaces $\Lambda_1$, $\Lambda_2$, while
$\Lambda_1$ and $\Lambda_2$ have the extra property that
they are orthogonal to each other. This means the inner
product of two arbitrary functions from $\Lambda_1$ and
$\Lambda_2$, respectively,
calculated as the covariance between them, is zero.
We then obtain its orthogonal complement
\[
\Lambda^\perp= \bigl\{\f(Y,\x)-E\bigl(\f\mid \bb^{\mathrm{T}}\x,Y
\bigr)\dvtx E(\f\mid\x)=E\bigl(\f\mid\bb^{\mathrm
{T}}\x\bigr), \forall \f \bigr
\}.
\]
The detailed derivation of $\Lambda$ and $\Lambda^\perp$ is
given in Appendix A.2 of~\cite{MaZhu2012}.
The form of $\Lambda^\perp$ permits many possibilities for
constructing
estimating equations.
For example, for
arbitrary functions $\g_i$ and $\ba_i$, the linear combination
\[
\sum_{i=1}^k \bigl\{\g_i
\bigl(Y,\bb^{\mathrm{T}}\x \bigr) -E \bigl(\g_i\mid\bb^{\mathrm{T}}
\x \bigr) \bigr\} \bigl\{\ba_i (\x ) - E \bigl(\ba_i\mid
\bb^{\mathrm
{T}}\x \bigr) \bigr\}
\]
will provide a
consistent semiparametric estimator since it is a valid
element in~$\Lambda^\perp$.
This form is exploited extensively in~\cite{MaZhu2012} to establish
links between the semiparametric approach and
various inverse regression methods.
Among all elements in $\Lambda^\perp$, the most interesting one is
the efficient score, defined as the orthogonal projection of
the score vector $\S_\beta$ onto $\Lambda^\perp$. We write the
efficient score as $\S\eff=\Pi(\S_\beta\mid\Lambda^\perp)$.
Because the
efficient score can be normalized to the efficient influence function,
it enables us to construct an efficient estimator of $\vecl(\bb)$ which
reaches the optimal semiparametric efficiency bound in the sense of
\cite{BickelKlaassenRitovWellner1993}.
In the supplementary document~\cite{supp}, we derive the efficient score function
to be
%
\begin{eqnarray}
\label{eqcseff}\qquad \S\eff\bigl(Y,\x,\bb^{\mathrm{T}}\x,\eta_2\bigr) =
\vecl \biggl[ \bigl\{\x-E\bigl(\x\mid \bb^{\mathrm{T}}\x\bigr)\bigr\}
\frac{\partial\log \{
\eta_2(Y,\bb^{\mathrm{T}}\x) \}}{\partial(\x^{\mathrm
{T}}\bb)} \biggr].
\end{eqnarray}
Hypothetically, the efficient estimator can be obtained through
implementing
\[
\sum_{i=1}^n\S\eff\bigl(Y_i,
\x_i,\bb^{\mathrm{T}}\x_i,\eta_2\bigr)=
\mathbf{0}.\vadjust{\goodbreak}
\]
However, $\S\eff$ is not readily implementable because it contains
the unknown quantities $E(\x\mid\x^{\mathrm{T}}\bb)$ and
$\partial\log\eta_2(Y, \bb^{\mathrm{T}}\x)/\partial(\x^{\mathrm
{T}}\bb)$.
For this reason, we first discuss a simpler alternative in the
following section.

\section{Locally efficient and efficient estimators}\label{seceff}
\subsection{Locally efficient estimators}\label{seceff1}
We now discuss how to construct a locally efficient estimator. This is
an estimator that contains some subjectively chosen components.
If the components are ``well'' chosen, the resulting estimator is
efficient. Otherwise, it is not efficient, but still consistent.
The efficient estimator defined in (\ref{eqcseff}) requires one to
estimate $\eta_2$,
the conditional p.d.f. of $Y$ on $\bb^{\mathrm{T}}\x$, and
its first derivative with respect to $\bb^{\mathrm{T}}\x$.
Although this is feasible, as we will describe in detail in Section
\ref{seceff2}, it certainly is not a trivial task as it involves several
nonparametric estimations.
Because of this, a compromise is to consider an
estimator that depends on a posited model of $\eta_2$.
Specifically, we would choose some favorite
form for $\eta_2$, denoted $\eta_2^*(Y,\bb^{\mathrm{T}}\x)$, and
utilize it
in place of $\eta_2$ to construct an estimating equation. If the posited
model is correct (i.e., $\eta_2^*=\eta_2$), then we would
have the optimal efficiency using the corresponding $\S\eff^*$.
However, even if the posited model is incorrect (i.e.,
$\eta_2^*\ne\eta_2$), we would still have consistency using the
corresponding $\S\eff^*$. A valid choice of $\S\eff^*$ that indeed
guarantees such property is\looseness=-1
\begin{eqnarray*}
&&\S\eff^*\bigl(Y_i,\x_i,\bb^{\mathrm{T}}
\x_i,\eta_2^*\bigr)\\
&&\qquad = \vecl \biggl( \bigl\{
\x_i- E\bigl(\x_i\mid \bb^{\mathrm{T}}\x_i
\bigr) \bigr\} \\
&&\hspace*{56pt}{}\times\biggl[ \frac{\partial\log \{\eta_2^*(Y_i,\bb^{\mathrm{T}}\x_i) \}
}{\partial(\x_i^{\mathrm{T}}\bb)} - E \biggl\{\frac{\partial\log
\eta_2^*(Y_i,\bb^{\mathrm{T}}\x_i)}{\partial(\x_i^{\mathrm
{T}}\bb)} \Big|\bb^{\mathrm{T}}
\x_i \biggr\} \biggr] \biggr).
\end{eqnarray*}\looseness=0
When $\eta_2^* = \eta_2$,
$E \{\partial\log
\eta_2^*(Y_i,\bb^{\mathrm{T}}\x_i)/\partial(\x_i^{\mathrm{T}}\bb
)\mid\bb^{\mathrm{T}}\x_i \}=\mathbf{
0}$, hence $\S\eff^*=\S\eff$. The
construction of a locally efficient estimator is often useful in
practice due
to its relative simplicity.
$\S\eff^*$ is almost readily applicable except that the two
expectations $E(\x_i\mid\bb^{\mathrm{T}}\x_i)$ and $E \{
\partial
\log
\eta_2^*(Y_i,\bb^{\mathrm{T}}\x_i)/\partial(\x_i^{\mathrm{T}}\bb
)\mid\bb^{\mathrm{T}}\x_i \}$
need to be estimated nonparametrically. One can use the familiar
kernel or local polynomial estimators. In Theorem
\ref{thlocal}, we show that under mild conditions, with the two
expectations estimated via the Nadaraya--Watson
kernel estimators, the local efficiency property indeed holds and
estimating the two expectations does not cause any difference from
knowing them in terms of its first order asymptotic
property.

We first
present the regularity conditions needed for the theoretical development.

\begin{longlist}

\item[(A1)] (\textit{The posited conditional density $\eta_2^*$}).
Denote $\u=\bb^{\mathrm{T}}\x$. The posited conditional density
$\eta_2^\ast(Y,\u)$ of $Y$ given $\u$ is bounded away from $0$ and
infinity on its support $\cal Y$. The second derivative of\vadjust{\goodbreak} $\log
\eta_2^\ast(Y,\u)$ with respect to $\u$
is continuous, positive definite and bounded.
In addition, there is an open
set $\bolds{\Omega}\in\mR^{p_t}$ which contains the true
parameter $\vecl(\bb)$, such that the third derivative of
$\eta_2(Y,\bb^{\mathrm{T}}\x)$ satisfies
\begin{eqnarray*}
\bigl\llvert \partial^3 \bigl\{\eta_2^\ast
\bigl(Y,\bb^{\mathrm{T}}\x\bigr) \bigr\} / \bigl(\partial \vecl(
\bb)_j\,\partial\vecl(\bb)_k\,\partial\vecl(
\bb)_l \bigr)\bigr\rrvert \le M_{jkl}^\ast(Y,\x)
\end{eqnarray*}
for all $\vecl(\bb)\in\bolds{\Omega}$ and $1\le j, k, l\le p_t$,
where $M_{jkl}^\ast(Y,\x)$ satisfies
$E \{{M_{jkl}^\ast}^2(Y,\break\x) \}<\infty$, and $\beta_j$ is the
$j$th component of
$\vecl(\bb)$.

\item[(A2)] (\textit{The nonparametric estimation}).
$E \{\partial\log
\eta_2^*(Y,\bb^{\mathrm{T}}\x)/\partial(\x^{\mathrm{T}}\bb)\mid
\bb^{\mathrm{T}}\x
\}$ and
$E(\x\mid\bb^{\mathrm{T}}\x)$
are estimated via the Nadaraya--Watson kernel estimator. For
simplicity, a
common
bandwidth $h$ is used which
satisfies $nh^{8}\to0$ and $nh^{2d}\to\infty$ as $n\to\infty$.

\item[(B1)](\textit{The true conditional density $\eta_2$}). The true
conditional density\break
$\eta_2(Y,\u)$ of $Y$ given $\u$ is bounded away from $0$ and
infinity on its support $\cal Y$. The first and second derivatives of
$\log\eta_2$ satisfy
\[
E \biggl[\frac{\partial \{\log\eta_2(Y,\bb^{\mathrm{T}}\x
) \}
}{\partial\vecl(\bb)} \biggr] = \mathbf{0}
\]
and
\begin{eqnarray*}
E \biggl[\frac{\partial \{\log\eta_2(Y,\bb^{\mathrm{T}}\x
) \}
}{\partial\vecl(\bb)} \frac{\partial \{\log\eta_2(Y,\bb^{\mathrm{T}}\x) \}
}{\partial
\vecl(\bb)^{\mathrm{T}}} \biggr] = -E \biggl[
\frac{\partial^2 \{\log\eta_2(Y,\bb^{\mathrm{T}}\x) \}} {
\partial\vecl(\bb)\,\partial\vecl(\bb)^{\mathrm{T}}} \biggr]
\end{eqnarray*}
is positive definite and bounded.
In addition, there is an open
set $\bolds{\Omega}\in\mR^{p_t}$ which contains the true
parameter $\vecl(\bb)$, such that the third derivative of
$\eta_2(Y,\bb^{\mathrm{T}}\x)$ satisfies
\begin{eqnarray*}
\bigl\llvert \partial^3 \bigl\{\eta_2\bigl(Y,
\bb^{\mathrm{T}}\x\bigr) \bigr\} / \bigl(\partial \vecl(\bb)_j
\,\partial\vecl(\bb)_k\,\partial\vecl(\bb)_l \bigr)\bigr
\rrvert \le M_{jkl}(Y,\x)
\end{eqnarray*}
for all $\vecl(\bb)\in\bolds{\Omega}$ and $1\le j, k, l\le p_t$,
where $M_{jkl}(Y,\x)$ satisfies
$E \{M_{jkl}^2(Y,\break \x) \}<\infty$, and $\beta_j$ is the
$j$th component of
$\vecl(\bb)$.

\item[(B2)](\textit{The bandwidths}). The bandwidths
satisfy $h_y\to0$, $b\to0$ and $h_x\to0$, and
$nh_y^{d+2}b\to\infty$,
$n^{1/2}\{h_x^2+(nh_x^d)^{-1/2}\}\{h_y^2+b^2+(nh_y^{d+2}b)^{-1/2}\}\to
0$.

\item[(C1)](\textit{The density functions of covariates}). Let $\u=
\bb^{\mathrm{T}}\x$. The density functions
of $\u$ and $\x$ are bounded away from $0$ and
infinity on their support $\cal U$ and $\cal X$ where ${\cal
U}= \{\u=\bb^{\mathrm{T}}\x\dvtx\x\in{\cal X} \}$ and
${\cal X}$ is
a compact support set of $\x$. Their second derivatives are finite
on their supports.

\item[(C2)](\textit{The smoothness}). The regression functions
$E(\x\mid\u)$ has a bounded and continuous derivative on $\cal U$.

\item[(C3)](\textit{The kernel function}). The univariate kernel function
$K(\cdot)$ is
a bounded symmetric probability density function, has a
bounded derivative and compact support $[-1,1]$, and satisfies $\mu_2
= \int u^2K(u)\,du\ne0$.
The $d$-dimensional kernel
function is a product of $d$ univariate kernel functions, that
is, $K(\u)=\prod_{j=1}^dK(u_j)$, and
$K_h(\u)=\prod_{j=1}^dK_h(u_j)=h^{-d}\prod_{j=1}^dK(u_j/h)$ for $\u
= (u_1,\ldots,u_d)^{\mathrm{T}}$ and any bandwidth $h$.
\end{longlist}

\begin{Th}\label{thlocal}
Under conditions \textup{(A1)--(A2)} and \textup{(C1)--(C3)}, the estimator obtained from
the estimating equation
\[
\sum_{i=1}^n \S\eff^*
\bigl(Y_i,\x_i,\bb^{\mathrm{T}}\x_i,
\eta_2^*,\widehat E\bigr)=\mathbf{0}
\]
is locally efficient. Specifically,
the estimator is consistent if $\eta_2^*\ne\eta_2$, and is
efficient if $\eta_2^*=\eta_2$. In addition, using the estimated
$\widehat
E(\cdot\mid\bb^{\mathrm{T}}\x)$ results in the same estimation
variance for
$\vecl(\bb)$ as using the true
$E(\cdot\mid\bb^{\mathrm{T}}\x)$. Specifically, the
estimate $\widehat\bb$ satisfies
\[
\sqrt{n}\bigl\{\vecl(\widehat{\bb})-\vecl(\bb)\bigr\}\to N\bigl\{\mathbf{0}, \A^{-1}\B\bigl(\A^{-1}\bigr)^{\mathrm{T}}\bigr\}
\]
when $n\to\infty$, where
\begin{eqnarray*}
\A=E \biggl\{\frac{\partial\S\eff^*(Y_i,\x_i,\bb^{\mathrm{T}}\x_i,\eta_2^*)} {
\partial\vecl(\bb)^{\mathrm{T}}} \biggr\},\qquad \B=E \bigl\{\S\eff^*
\bigl(Y_i,\x_i,\bb^{\mathrm{T}}\x_i,
\eta_2^*\bigr)^{\otimes2} \bigr\}.
\end{eqnarray*}
\end{Th}
In Theorem~\ref{thlocal} and thereafter, we use $\mathbf{v}^{\otimes2}$
to denote
$\mathbf{v}\mathbf{v}^{\mathrm{T}}$ for any matrix or
vector~$\mathbf{v}$,
and use
$\widehat E$ to denote the nonparametrically estimated expectation.

We describe how to implement the locally efficient estimator in several
specific cases.
For example, when $Y$ is continuous, we can propose a simple
conditional normal model for $\eta_2$
and hence obtain the locally efficient estimator based on summing terms
of the form
%
\begin{eqnarray}\label{eqlocal}
&&\n \S\eff^*\bigl(Y,\x,\bb^{\mathrm{T}}\x,\eta_2^*\bigr)
\nonumber
\\[-8pt]
\\[-8pt]
\nonumber
&&\qquad= \vecl \biggl( \bigl\{\x- E\bigl(\x\mid \bb^{\mathrm{T}}\x
\bigr) \bigr\} \biggl[\bigl\{Y-E\bigl(Y\mid\bb^{\mathrm{T}}\x \bigr)\bigr\}
\frac
{\partial
E^*(Y\mid\bb^{\mathrm{T}}\x)}{\partial (\x^{\mathrm{T}}\bb
)} \biggr] \biggr)
\end{eqnarray}
evaluated at different observations. Here
$E^*(\cdot\mid\bb^{\mathrm{T}}\x)$ is computed using the model~$\eta_2^*$.
When $Y$ is binary, a common model to posit for $\eta_2$ is a logistic
model. The summation of the terms of form (\ref{eqlocal})
evaluated at different observations also provides a locally efficient
estimator.
When $Y$ is a counting response variable,
the Poisson model is a popular choice for $\eta_2$.
This choice also yields an identical locally efficient estimator formed by
the sum of (\ref{eqlocal}).
The benefits of these locally efficient estimators are two-fold.
The first benefit lies in the robustness property, in
that they guarantee the consistency of the resulting estimators
regardless of the proposed model. The second benefit is
their computational simplicity
gained through avoiding estimating the conditional density $\eta_2$
and its derivative.
In addition, if, by luck, the posited model happens to be correct, then
the estimator is efficient.

\begin{Rem}
We have restricted the posited model $\eta_2^*$ to be a
completely known model in order to illustrate the local efficiency
concept. In fact, one can also posit\vadjust{\goodbreak} a model $\eta_2^*$ that contains
an additional unknown parameter vector, say~$\bg$. As long as $\bg$
can be
estimated at the root-$n$ rate, the resulting estimator with the estimator
$\widehat\bg$ plugged in is also referred to as a locally efficient
estimator. In addition, if model $\eta_2^*$ contains the true~$\eta_2$, say
$\eta_2^*(Y, \bb^{\mathrm{T}}\x, \bg_0)=\eta_2(Y, \bb^{\mathrm
{T}}\x)$, and
$\bg_0$
is estimated consistently by $\widehat\bg$ at the root-$n$ rate, then
the resulting estimator $\S\eff^*$ with $\eta_2^*(Y,\bb^{\mathrm
{T}}\x
,\widehat\bg)$ plugged in
is efficient.
\end{Rem}

\begin{Rem}
Even if efficiency is not sought after and consistency is the sole
purpose, at least one nonparametric operation, such as one that relates to
estimating $E(\x\mid\bb^{\mathrm{T}}\x)$, is needed. Thus, to
completely avoid
nonparametric procedures, the only option is to impose additional
assumptions. The most popular linearity condition in the literature
assumes
$E(\x\mid\bb^{\mathrm{T}}\x)=\bb(\bb^{\mathrm{T}}\bb)^{-1}\bb^{\mathrm{T}}\x$.
Since Theorem
\ref{thlocal} allows an arbitrary $\eta^*$, the most obvious choice
in practice is probably the exponential link
functions. For example, if we choose $\eta_2^*$ to be the normal
link function when $d=1$, then the locally efficient estimator
degenerates to a simple form, where
\[
\S\eff^*=\vecl \bigl[\bigl\{\x-\bb\bigl(\bb^{\mathrm{T}}\bb\bigr)^{-1}
\bb^{\mathrm{T}}\x\bigr\} \bigl(Y-\bb^{\mathrm{T}}\x\bigr) \bigr].
\]
If we are even bolder and decide to replace $Y-\bb^{\mathrm{T}}\x$ with
$Y$, which is still valid given that the first term alone already guarantees
consistency under the linearity condition, then we obtain the ordinary
least square estimator~\cite{LiDuan1991}. Further connections to
other existing methods are
elaborated in~\cite{MaZhu2012}.
\end{Rem}

\subsection{The efficient estimator}\label{seceff2}
Now we pursue the truly efficient
estimator that reaches the semiparametric efficiency bound. This is
important because
in terms of reaching the optimal efficiency, relying on a
posited model $\eta_2^*$ to be true or to contain the true $\eta_2$ is
not a satisfying practice. Intuitively, it is
easy to imagine that in constructing the locally efficient estimator,
if we posit a larger model $\eta_2^*$, the chance of it containing the
true model $\eta_2$ becomes larger, hence the chance of reaching the
optimal efficiency also increases. Thus, if we can propose the
``largest'' possible model for $\eta_2^*$, we will guarantee to have
$\eta_2^*$ containing
$\eta_2$. If we can also estimate the parameters in $\eta_2^*$
``correctly,'' we will then guarantee the efficiency. This ``largest''
model with a ``correctly'' estimated parameter turns
out to be what the nonparametric estimation is able to provide.
This amounts to estimating $E(\x\mid\bb^{\mathrm{T}}\x)$, $\eta_2$ and its
first derivative nonparametrically
in (\ref{eqcseff}).

We first discuss how to estimate $\eta_2$ and its first derivative,
based on $(Y_i, \bb^{\mathrm{T}}\x_i),  i=1, \ldots, n$.
This is a problem of estimating conditional density and its derivative.
We use
the idea of the ``double-kernel'' local linear smoothing method studied
in~\cite{FanYaoTong1996}. Consider $K_b(Y-y) =
b^{-1}K \{(Y-y)/b \}$ with $y$ running through all possible
values, where $K(\cdot)$ is a symmetric density function, and $b>0$ is
a bandwidth. Then
$E \{K_b(Y-y)\mid\bb^{\mathrm{T}}\x \}$ converges to
$\eta_2(y,\bb^{\mathrm{T}}\x)$ as
$b$ tends to 0. This observation motivates us
to estimate $\eta_2$ and its first derivative, evaluated at $(y,
\bb^{\mathrm{T}}\x)$ through minimizing the
following weighted least squares:
\begin{eqnarray*}
\sum_{i=1}^n \bigl\{K_b(Y_i-y)-a-
\mathbf{b}^{\mathrm{T}}\bigl(\bb^{\mathrm{T}}\x_i-\bb^{\mathrm{T}}
\x \bigr) \bigr\}^2 K_{h_y}\bigl(\bb^{\mathrm{T}}
\x_i-\bb^{\mathrm{T}}\x\bigr),
\end{eqnarray*}
where $h_y$ is a bandwidth, and $K_{h_y}$ is a multivariate kernel
function.
The minimizers~$\widehat a$ and $\widehat{\mathbf{b}}$ are the estimators
of $\eta_2$
and $\partial\eta_2/\partial(\bb^{\mathrm{T}}\x)$.
Let the resulting estimators be $\widehat\eta_2(\cdot)$ and
$\widehat\eta_2'(\cdot)$.

It remains to estimate $E(\x\mid\bb^{\mathrm{T}}\x)$. Using the
Nadaraya--Watson
kernel estimator, we have
\[
\widehat E\bigl(\x\mid\bb^{\mathrm{T}}\x\bigr) =\frac{\sum_{i=1}^n\x_i
K_{h_x}(\bb^{\mathrm{T}}\x_i-\bb^{\mathrm{T}}\x)} {
\sum_{i=1}^n
K_{h_x}(\bb^{\mathrm{T}}\x_i-\bb^{\mathrm{T}}\x)},
\]
where $h_x$ is a bandwidth, and $K_{h_x}$ is a multivariate kernel
function.
The algorithm for obtaining the efficient estimator is the following:
\begin{itemize}
\item \textit{Step} 1.
Obtain an initial root-$n$ consistent
estimator of $\bb$, denoted as $\widetilde\bb$, through, for example,
a simple
locally efficient estimation procedure from Section~\ref{seceff1}.
\item \textit{Step} 2.
Perform nonparametric estimation of $\eta_2(Y,\widetilde\bb^{\mathrm
{T}}\x
)$ and
its first derivative $\partial
\{\eta_2(Y,\widetilde\bb^{\mathrm{T}}\x)\}/\partial(\widetilde
\bb^{\mathrm{T}}\x)$.
Write the resulting estimators as $\widehat\eta_2(\cdot)$ and
$\widehat\eta_2'(\cdot)$.
\item \textit{Step} 3.
Perform nonparametric estimation of $E(\x\mid\widetilde\bb^{\mathrm
{T}}\x)$.
Write the resulting estimator as $\widehat E(\cdot)$.
\item \textit{Step} 4.
Plug $\widehat\eta_2(Y, \bb^{\mathrm{T}}\x)$, $\widehat\eta_2'(Y,\bb^{\mathrm{T}}\x)$ and $\widehat
E(\x\mid\bb^{\mathrm{T}}\x)$ into $\S\eff$ and solve the
estimating equation
\[
\sum_{i=1}^n\S\eff\bigl(Y_i,
\x_i,\bb^{\mathrm{T}}\x_i,\widehat\eta_2,
\widehat \eta_2',\widehat E\bigr)=\mathbf{0}
\]
to obtain the efficient
estimator $\widehat\bb$.
\end{itemize}

In performing the various nonparametric estimations in
steps 2 and 3, as well as in obtaining the locally efficient
estimator in Section~\ref{seceff1}, bandwidths need to be
selected. Because the final estimator is very insensitive to the
bandwidths, as indicated by conditions (A2), (B2)
and Theorems~\ref{thlocal},~\ref{theff}, where a range of different
bandwidths all lead to the same asymptotic property of the final
estimator, we suggest that one should select the corresponding bandwidths
by taking the sample size $n$ to its suitable power to satisfy (B2),
and then multiply a constant to scale it, instead of
performing a full-scale cross validation procedure.
For example, when $d=1$, we let $h=n^{-1/5}, h_x=n^{-1/5},
h_y=n^{-1/6}, b=n^{-1/7}$, and when $d=2$, we let $h=n^{-1/6}, h_x=n^{-1/6},
h_y=n^{-1/7}, b=n^{-1/8}$, each multiplied by the standard deviation of
the regressors calculated at the current $\widehat\bb$ value.

The estimator from the above algorithm, $\widehat\bb$, with
its upper $d\times d$ submatrix being~$\I_d$,
reaches the optimal semiparametric efficiency bound. We present this
result in
Theorem~\ref{theff}.

\begin{Th}\label{theff}
Under conditions \textup{(B1)--(B2)} and \textup{(C1)--(C3)}, the estimator obtained from
the estimating equation
\[
\sum_{i=1}^n \S\eff\bigl(Y_i,
\x_i,\bb^{\mathrm{T}}\x_i,\widehat\eta_2,
\widehat \eta_2',\widehat E\bigr)=\mathbf{0}
\]
is efficient. Specifically, when $n\to\infty$, the estimator
of $\vecl(\bb)$ satisfies
\begin{eqnarray*}
\sqrt{n}\bigl\{{\vecl(\widehat\bb)}-\vecl(\bb)\bigr\}\to N \bigl(\mathbf{0}, \bigl[E
\bigl\{\S\eff\bigl(Y,\x,\bb^{\mathrm{T}}\x,\eta_2
\bigr)^{\otimes
2} \bigr\} \bigr]^{-1} \bigr)
\end{eqnarray*}
in distribution.
\end{Th}

\begin{Rem}
It is discovered that for
certain p.d.f. $\eta_2$, such as when the inverse mean function $E(\x
\mid Y)$ degenerates,
some inverse, regression-based methods, such as SIR, would fail to
exhaustively recover $\cs$. However, this is not the case for the
efficient estimator
proposed here. That is, our proposed efficient estimator,
similar to dMAVE~\cite{Xia2007},
has the exhaustiveness property~\cite{LiZhaChiaromonte2005}.
In fact, as it is listed in the regularity conditions, as
long as
the asymptotic covariance matrix is not
singular and is bounded away from infinity, our method is always able
to produce the efficient estimator.
\end{Rem}

\begin{Rem}
It can be easily verified that the above efficient asymptotic
variance-covariance matrix can be explicitly written out as
\begin{eqnarray*}
&&E \bigl\{\S\eff\bigl(Y,\x,\bb^{\mathrm{T}}\x,\eta_2
\bigr)^{\otimes
2} \bigr\}
\\
&&\qquad = E \biggl(E \biggl[ \biggl\{\frac{\partial\log\eta_2(Y,
\bb^{\mathrm{T}}\x)}{\partial(\bb^{\mathrm{T}}\x)} \biggr\}^{\otimes2}\Big\mid
\bb^{\mathrm{T}}\x \biggr] \otimes E \bigl[ \bigl\{\x_{l}-E\bigl(
\x_{l}\mid \bb^{\mathrm{T}}\x\bigr) \bigr\}^{\otimes2}\mid
\bb^{\mathrm{T}}\x \bigr] \biggr),
\end{eqnarray*}
where 
$\x_l$ is the vector
formed by the lower
$p-d$ components of $\x$.
Thus, the asymptotic variance of $\vecl(\widehat\bb)$
is nonsingular as long as
both $E [ \{\partial\log\eta_2(Y,\break
\bb^{\mathrm{T}}\x)/ {\partial(\bb^{\mathrm{T}}\x)} \}^{\otimes2}\mid
\bb^{\mathrm{T}}\x ]$ and $E [ \{\x_{l}-E(\x_{l}\mid
\bb^{\mathrm{T}}\x) \}^{\otimes2}\mid\bb^{\mathrm{T}}\x
]$ are
nonsingular. The nonsingularity of the first matrix is a standard
requirement on the information matrix of the true model $\eta_2$ and
is usually satisfied. On the other hand, $E (E [ \{\x_{l}-E(\x_{l}\mid
\bb^{\mathrm{T}}\x) \}^{\otimes2}\mid\bb^{\mathrm{T}}\x
] )$
is always
guaranteed to be nonsingular. This is because if it is singular, then
there exists a unit vector $\ba$ with the first $d$ components zero,
such that $\ba^{\mathrm{T}}\x$ is a deterministic
function of $\bb^{\mathrm{T}}\x$.
This violates our assumption that $\ba^{\mathrm{T}}\x$ cannot be a
deterministic function of $\bb^{\mathrm{T}}\x$
unless $\ba$ lies within the column space of $\bb$.
\end{Rem}
%

\section{Simulation study}\label{secsimu}

In this section we conduct simulations to evaluate the finite sample performance
of our efficient and locally efficient estimators and compare them with
several existing methods.

We consider the following three examples:
\begin{longlist}[(1)]
\item[(1)] We generate $Y$ from a normal population with mean function
$\x^{\mathrm{T}}\bb$ and variance $1$.
\item[(2)] We generate $Y$ from a normal population with mean
function\break
$\sin (2\x^{\mathrm{T}}\bb )+2\exp (2+\x^{\mathrm
{T}}\bb
)$
and variance function
$\log\{2+(\x^{\mathrm{T}}\bb)^2\}$.
\item[(3)] We generate $Y$ from a normal population with mean function
$2 (\x^{\mathrm{T}}\bb_1 )^{2}$ and variance function
$2\exp
(\x^{\mathrm{T}}\bb_2 )$.
\end{longlist}

In the simulated examples 1 and 2, we set
$\bb=(1.3, -1.3, 1.0, -0.5, 0.5,\break -0.5)^{\mathrm{T}}$ and
generate $\x= (X_1,\ldots,X_6)^{\mathrm{T}}$
as follows. We generate
$X_1$, $X_2$, $e_1$ and $e_2$ independently from a
standard normal distribution, and form $X_3=0.2X_1+0.2(X_2+2)^2+0.2e_1$,
$X_4=0.1+0.1(X_1+X_2)+0.3(X_1+1.5)^2+0.2e_2$. We generate $X_5$ and
$X_6$ independently from
Bernoulli distributions with success probability $\exp (X_1
)/ \{1+
\exp (X_1 ) \}$ and
$\exp (X_2 )/ \{1+
\exp (X_2 ) \}$, respectively.

Example 3 follows the setup of
Example 4.2 in~\cite{Xia2007}.
In this example, we set
$\bb_1 = (1, 2/3, 2/3, 0, -1/3, 2/3)^{\mathrm{T}}$
and $\bb_2 = (0.8,0.8,-0.3,0.3,0,0)^{\mathrm{T}}$.
We form the
covariates $\x$ by setting $X_1=U_1-U_2$, $X_2=U_2-U_3-U_4$,
$X_3=U_3+U_4$, $X_4=2U_4$, $X_5=U_5+0.5U_6$ and $X_6=U_6$, where
$U_1$ is generated from a Bernoulli distribution with probability
0.5 to be 1 or $-1$, $U_2$ is also generated from Bernoulli
distribution, with probability 0.7 to be $\sqrt{3/7}$ and probability
$0.3$ to be $-\sqrt{7/3}$. The remaining four components of $\u$ are generated
from a uniform distribution between $-\sqrt{3}$ and $\sqrt{3}$.
The six components of $\u=(U_1,\ldots,U_6)^{\mathrm{T}}$ are independent,
marginally having zero mean and unit variance.
We construct $\x$ through $\u$ in this way to allow the components of
$\x$ to be correlated.

For the purpose of comparison, we implement six estimators:
``Oracle,'' ``Eff,'' ``Local,'' ``dMAVE,'' ``SIR'' and ``DR.''
The names of the estimators suggest
the nature of these estimators, while we briefly explain them in the following:
\begin{description}
\item Oracle: the oracle estimate which correctly specifies
$\eta_2$ in (\ref{eqcseff}), but we estimate $E(\x\mid\bb^{\mathrm{T}}
\x)$ through
kernel regressions. We remark here that the oracle estimator is
not a realistic estimator because
$\eta_2$ is usually unknown. We include the oracle estimator here
to provide a benchmark since this is the best performance one could
hope for.

\item Eff: the efficient estimator which estimates $E(\x\mid\bb^{\mathrm{T}}\x)$,
$\eta_2$ and $\eta_2'$ through
nonparametric regressions. See Section~\ref{seceff2} for a description
about this efficient estimator.

\item Local: the locally efficient estimate which mis-specifies
the model $\eta_2$, and estimates $E(\cdot\mid\bb^{\mathrm{T}}\x
)$ through
nonparametric regression. This is an implementation of (\ref{eqlocal}).

\item dMAVE: the conditional density based
minimum average variance estimation proposed by~\cite{Xia2007}.

\item SIR: the sliced inverse regression~\cite{Li1991} which
estimates $\bb$ as the first $d$ principal eigenvectors
of $\bSig^{-1}\cov \{E (\x\mid Y ) \}\bSig^{-1}$,
where $\bSig= \cov(\x)$.

\item DR: the directional regression~\cite{LiWang2007} which
estimates $\bb$ as the first $d$ principal eigenvectors
of the kernel matrix $\bSig^{-1/2}
E\{2\I_p-\A(Y,\tY)\}^2\bSig^{-1/2}$, where
$\A(Y,\tY)=\bSig^{-1/2}E\{(\x-\tx)(\x-\tx)^{\mathrm
{T}}\mid
Y,\tY\}\bSig^{-1/2}$, and $(\tx,\tY)$ is an independent
copy of $(\x,Y)$.
\end{description}

%
\begin{table}
\caption{The average (``\textup{ave}'') and
the sample standard errors (``\textup{std}'') for various estimates,
and the inference results, respectively, the average of the estimated
standard deviation
(``$\widehat{\mathrm{std}}$'') and the~coverage of the estimated 95\%
confidence
interval (``95\%''),
of the oracle estimator and the efficient estimator, of
$\bb$ in simulated example~1}\label{tabletab1}
\begin{tabular*}{\textwidth}{@{\extracolsep{\fill
}}lcd{1.4}d{2.4}d{1.4}d{2.4}d{1.4}d{2.4}@{}}
\hline
&&\multicolumn{1}{c}{$\bolds{\beta_1}$}& \multicolumn{1}{c}{$\bolds
{\beta_2}$}& \multicolumn{1}{c}{$\bolds{\beta_3}$}
&\multicolumn{1}{c}{$\bolds{\beta_4}$}&
\multicolumn{1}{c}{$\bolds{\beta_5}$}&\multicolumn{1}{c@{}}{$\bolds
{\beta_6}$}\\
&&\multicolumn{1}{c}{\textbf{1.3}} & \multicolumn{1}{c}{$\bolds
{-1.3}$} & \multicolumn{1}{c}{\textbf{1}} &
\multicolumn{1}{c}{$\bolds{-0.5}$} & \multicolumn{1}{c}{\textbf
{0.5}} &\multicolumn{1}{c@{}}{$\bolds{-0.5}$} \\
\hline
Oracle&ave&1.2978&-1.3036&1.0049&-0.4985&0.5033&-0.4943\\
&std&0.1221&0.1477&0.1505&0.1169&0.0966&0.1049\\
&$\widehat{\mathrm{std}}$&0.1264&0.1510&0.1527&0.1212&0.0983&0.1052\\
&95\%&0.9510&0.9540&0.9440&0.9540&0.9520&0.9450\\[3pt]
Eff&ave&1.2980&-1.3046&1.0064&-0.4990&0.5040&-0.4936\\
&std&0.1280&0.1546&0.1567&0.1221&0.1000&0.1075\\
&$\widehat{\mathrm{std}}$&0.1317&0.1588&0.1602&0.1264&0.1011&0.1084\\
&95\%&0.9480&0.9380&0.9380&0.9440&0.9480&0.9510\\[3pt]
Local&ave&1.3052 & -1.2629 & 0.9687 & -0.4988 & 0.5023 & -0.4897 \\
&std&0.1478 & 0.1736 & 0.1715 & 0.1393 & 0.1069 & 0.1153 \\[3pt]
dMAVE&ave&1.2599 & -1.2933 & 1.0014 & -0.4763 & 0.4984 & -0.4935 \\
&std&0.1932 & 0.1427 & 0.1550 & 0.1701 & 0.1368 & 0.1378 \\[3pt]
SIR&ave&1.3881 & -1.1930 & 0.9261 & -0.5968 & 0.4793 & -0.4724 \\
&std&0.1696 & 0.1522 & 0.1414 & 0.1489 & 0.0976 & 0.0995 \\[3pt]
DR&ave&0.9935 & -0.2217 & 0.1930 & -0.6863 & 0.1245 & -0.1071 \\
&std&0.6567 & 1.2305 & 1.0107 & 0.6411 & 0.3069 & 0.2999 \\
\hline
\end{tabular*}
\end{table}

We repeat each experiment 1000 times with sample size
$n=500$. The results are summarized
in Table~\ref{tabletab1} for example 1,
Table~\ref{tabletab2} for example 2 and Table
\ref{tabletab3}
for example 3.
Because the estimators we propose here use a different
parameterization of the central subspace $\cs$ from the existing
methods such as SIR, DR or dMAVE, we transform the results from all the
estimation procedures to
the original $\bb$ used to generate the data
for a fair and intuitive comparison.

%
\begin{table}
\caption{The average (``\textup{ave}'') and
the sample standard errors (``\textup{std}'') for various estimates,
and the inference results, respectively, the average of the estimated
standard deviation
(``$\widehat{\mathrm{std}}$'') and the coverage of the estimated 95\%
confidence
interval (``95\%''),
of the oracle estimator and the~efficient estimator, of
$\bb$ in simulated example~2}\label{tabletab2}
\begin{tabular*}{\textwidth}{@{\extracolsep{\fill
}}lcd{1.4}d{2.4}d{1.4}d{2.4}d{1.4}d{2.4}@{}}
\hline
&& \multicolumn{1}{c}{$\bolds{\beta_1}$}& \multicolumn
{1}{c}{$\bolds{\beta_2}$}&
\multicolumn{1}{c}{$\bolds{\beta_3}$}&\multicolumn{1}{c}{$\bolds
{\beta_4}$}&\multicolumn{1}{c}{$\bolds{\beta_5}$}&
\multicolumn{1}{c@{}}{$\bolds{\beta_6}$}\\
&&\multicolumn{1}{c}{\textbf{1.3}} & \multicolumn{1}{c}{$\bolds
{-1.3}$} & \multicolumn{1}{c}{\textbf{1}} &
\multicolumn{1}{c}{$\bolds{-0.5}$} & \multicolumn{1}{c}{\textbf
{0.5}} &
\multicolumn{1}{c@{}}{$\bolds{-0.5}$} \\
\hline
Oracle&ave&1.2999&-1.3001&1.0001&-0.4999&0.5002&-0.4999\\
&std&0.0023&0.0025&0.0028&0.0022&0.0023&0.0024\\
&$\widehat{\mathrm{std}}$&0.0021&0.0020&0.0026&0.0020&0.0021&0.0023\\
&95\%&0.9260&0.9070&0.9270&0.9220&0.9210&0.9380\\[3pt]
Eff&ave&1.2996&-1.2999&0.9998&-0.4996&0.5002&-0.5000\\
&std&0.0116&0.0116&0.0117&0.0111&0.0068&0.0079\\
&$\widehat{\mathrm{std}}$&0.0123&0.0124&0.0124&0.0120&0.0075&0.0081\\
&95\%&0.9480&0.9550&0.9570&0.9450&0.9630&0.9520\\[3pt]
Local&ave&1.2992 & -1.3010 & 1.0007 & -0.4993 & 0.5011 & -0.5001 \\
&std&0.0155 & 0.0210 & 0.0209 & 0.0140 & 0.0142 & 0.0147 \\[3pt]
dMAVE&ave&1.2405 & -1.3422 & 1.0303 & -0.4490 & 0.5114 & -0.5134 \\
&std&0.0229 & 0.0151 & 0.0133 & 0.0153 & 0.0081 & 0.0082 \\[3pt]
SIR&ave&0.3064 & -1.6387 & 1.2390 & 0.2477 & 0.4697 & -0.4743 \\
&std&0.1248 & 0.3965 & 0.3149 & 0.1057 & 0.1135 & 0.1141 \\[3pt]
DR&ave&0.3424 & 0.8686 & -0.6620 & -0.6895 & -0.1923 & 0.1912 \\
&std&0.2550 & 1.2518 & 0.9653 & 0.6938 & 0.3360 & 0.3410 \\
\hline
\end{tabular*}
\end{table}

\begin{sidewaystable}
\tablewidth=\textheight
\tablewidth=\textwidth
\caption{The average (``\textup{ave}'') and
the sample standard errors (``\textup{std}'') for various estimates,
and the inference results, respectively, the average of the estimated
standard deviation
(``$\widehat{\mathrm{std}}$'') and the coverage of the estimated 95\%
confidence
interval (``95\%''),
of the oracle estimator and the efficient estimator, of
$\bb$ in simulated example~3}\label{tabletab3}
\begin{tabular*}{\textwidth}{@{\extracolsep{\fill
}}lcd{1.4}d{1.4}d{2.4}d{2.4}d{2.4}d{1.4}d{1.4}d{1.4}d{2.4}d{1.4}d{2.4}d{2.4}@{}}
\hline
&&\multicolumn{1}{c}{$\bolds{\beta_{11}}$}&\multicolumn
{1}{c}{$\bolds{\beta_{21}}$}&\multicolumn{1}{c}{$\bolds{\beta_{31}}$}&
\multicolumn{1}{c}{$\bolds{\beta_{41}}$}&\multicolumn{1}{c}{$\bolds
{\beta_{51}}$}&\multicolumn{1}{c}{$\bolds{\beta_{61}}$}&\multicolumn
{1}{c}{$\bolds{\beta_{12}}$}&\multicolumn{1}{c}{$\bolds{\beta_{22}}$}&
\multicolumn{1}{c}{$\bolds{\beta_{32}}$}&\multicolumn{1}{c}{$\bolds
{\beta_{42}}$}&\multicolumn{1}{c}{$\bolds{\beta_{52}}$}&\multicolumn
{1}{c@{}}{$\bolds{\beta_{62}}$}\\
&&\multicolumn{1}{c}{\textbf{1}} &\multicolumn{1}{c}{\textbf
{0.6667}}&\multicolumn{1}{c}{\textbf{0.6667}}&\multicolumn
{1}{c}{\textbf{0}}&
\multicolumn{1}{c}{$\bolds{-0.3333}$}&\multicolumn{1}{c}{\textbf
{0.6667}}&\multicolumn{1}{c}{\textbf{0.8}}&\multicolumn
{1}{c}{\textbf{0.8}}&
\multicolumn{1}{c}{$\bolds{-0.3}$}&\multicolumn{1}{c}{\textbf
{0.3}}&\multicolumn{1}{c}{\textbf{0}}&\multicolumn{1}{c@{}}{\textbf
{0}}\\
\hline
Oracle &
ave&1.0009&0.6676&0.6674&0.0002&-0.3339&0.6675&0.8064&0.8064&-0.2905&0.2969&-0.0047&0.0053\\
&
std&0.0305&0.0305&0.0325&0.0099&0.0198&0.0314&0.0860&0.0860&0.0902&0.0291&0.0550&0.0854\\
& $\widehat{\mathrm
{std}}$&0.0275&0.0275&0.0295&0.0109&0.0178&0.0276&0.0828&0.0828&0.0876&0.0296&0.0547&0.0826\\
& 95\%
&0.9270&0.9270&0.9300&0.9590&0.9200&0.9110&0.9410&0.9410&0.9320&0.9450&0.9520&0.9430\\
[3pt]
Eff
&ave&1.0097&0.6763&0.6764&-0.0000&-0.3384&0.6752&0.8038&0.8038&-0.3067&0.3105&0.0022&-0.0003\\
&std&0.0714&0.0714&0.0745&0.0162&0.0434&0.0740&0.1737&0.1737&0.1993&0.0485&0.1511&0.1895\\
&$\widehat{\mathrm
{std}}$&0.0709&0.0709&0.0734&0.0175&0.0454&0.0702&0.1439&0.1439&0.1490&0.0381&0.0973&0.1439\\
&95\%
&0.9280&0.9280&0.9350&0.9530&0.9460&0.9430&0.9230&0.9230&0.9240&0.9410&0.9150&0.9080\\
[3pt]
local&ave&1.0633&0.7300&0.7372&-0.0072&-0.3701&0.7468&0.7689&0.7689&-0.3066&0.2754&-0.0116&-0.0042\\
&std&1.8783&1.8783&2.1273&0.2493&1.0694&2.3913&1.1281&1.1281&1.5767&0.4517&0.2192&0.2516\\
[3pt]
dMAVE
&ave&0.8884&0.6079&-0.1703&0.2119&-0.2498&0.5065&0.8282&0.7722&-0.0901&0.2371&-0.0153&0.0354\\
&std&0.0748&0.1021&0.0951&0.0569&0.0888&0.1155&0.0379&0.0378&0.1188&0.0731&0.0761&0.0489\\
[3pt]
SIR
&ave&0.5443&0.3781&-0.3301&0.1816&-0.0944&0.1976&0.7768&0.6849&-0.4083&0.2908&0.0441&-0.0828\\
&std&0.1514&0.1414&0.0863&0.0586&0.1257&0.2022&0.0650&0.0808&0.1098&0.0748&0.1059&0.0831\\
[3pt]
DR
&ave&0.6332&0.2753&-0.2968&0.0939&-0.2701&0.5422&0.7004&0.6823&-0.4512&0.1498&0.0013&-0.0151\\
&std&0.1813&0.2009&0.1003&0.0739&0.1288&0.1567&0.1063&0.1446&0.1688&0.0880&0.1639&0.0945\\
\hline
\end{tabular*}
\end{sidewaystable}

From the results in Table~\ref{tabletab1},
we can see that Oracle, Eff, Local, dMAVE provide estimators with
small bias, while
SIR and DR have substantial bias in some of the elements in
$\bb$. For example, the average of the second estimated component of
$\bb$
obtained by
DR is $-0.2217$, in contrast to the true value $-1.3$.
This is because the covariate $\x$ does not satisfy the linearity
or the constant variance condition, and hence violates the requirement
of SIR and DR.
Although Local and dMAVE both appear consistent, they have much
larger variance in some components than Eff. For example, in estimating
$\beta_1$,
the asymptotic variance of dMAVE is 0.1932, whereas that of Eff is as
small as 0.1264.
This is not surprising
since Eff is asymptotically efficient. In fact, for this very simple
setting, the estimation variance of Eff is almost as good as Oracle,
which indicates that the asymptotic efficiency already exhibits for
$n=500$.

We also provide the average of the estimated standard error using
the results in Theorem~\ref{theff} and the 95\% coverage in
Table~\ref{tabletab1}. The numbers show a close approximation of the
sample and estimated standard error and 95\% coverage is reasonable close
to the nominal value.

Similar phenomena are observed for the simulated example 2 from Table~\ref{tabletab2},
where SIR and DR are
biased, Local and dMAVE are consistent but have larger
variability than Eff and Oracle. In this more complex model where the
mean function is
highly nonlinear and the error is heteroscedastic, we lose the
proximity between the oracle performance and the Eff
performance. This is probably because $n=500$ is still too small for
this model. The inference
results in Table~\ref{tabletab2}, however, are still satisfactory, indicating
that although we cannot achieve the theoretical optimality, inference
is still sufficiently reliable.

What we observe in Table~\ref{tabletab3}, for the simulated example
3, tells a
completely different story. For this case with $d=2$, both the
linearity and
the constant variance condition are violated. In addition,
$\x$ contains categorical variables. dMAVE, SIR and DR
all fail to provide good estimators in terms of estimation bias. Local
and Eff remain to be
consistent, although like in the simulated example 2, we can no longer
hope to
see the optimality as the estimation standard error is much larger
than the Oracle estimator. Inference results presented in Table
\ref{tabletab3} still show satisfactory 95\% coverage values, while
the average estimated estimation standard error can deviate away from
the sample standard error. This is caused by some numerical
instability of a small proportion of the simulation repetitions. In fact,
if we replace the average with the median estimated standard error,
the results are closer.

\section{An application}\label{secexample}

We use the proposed efficient estimator
to analyze a dataset concerning the employees' salary in
the Fifth National Bank of Springfield~\cite
{AlbrightWinstonZappe1999}. The
aim of the study is to understand how an employee's
salary associates with his/her social characteristics. We regard an
employee's annual salary as the response variable $Y$, and several
social characteristics as the associated covariates. These
covariates are, specifically, current job level
($X_1$); number of years working at the bank ($X_2$);
age ($X_3$); number of years working at other banks ($X_4$);
gender ($X_5$); whether the job is computer
related ($X_6$). After removing an obvious outlier, the dataset
contains 207 observations.

\begin{figure}

\includegraphics{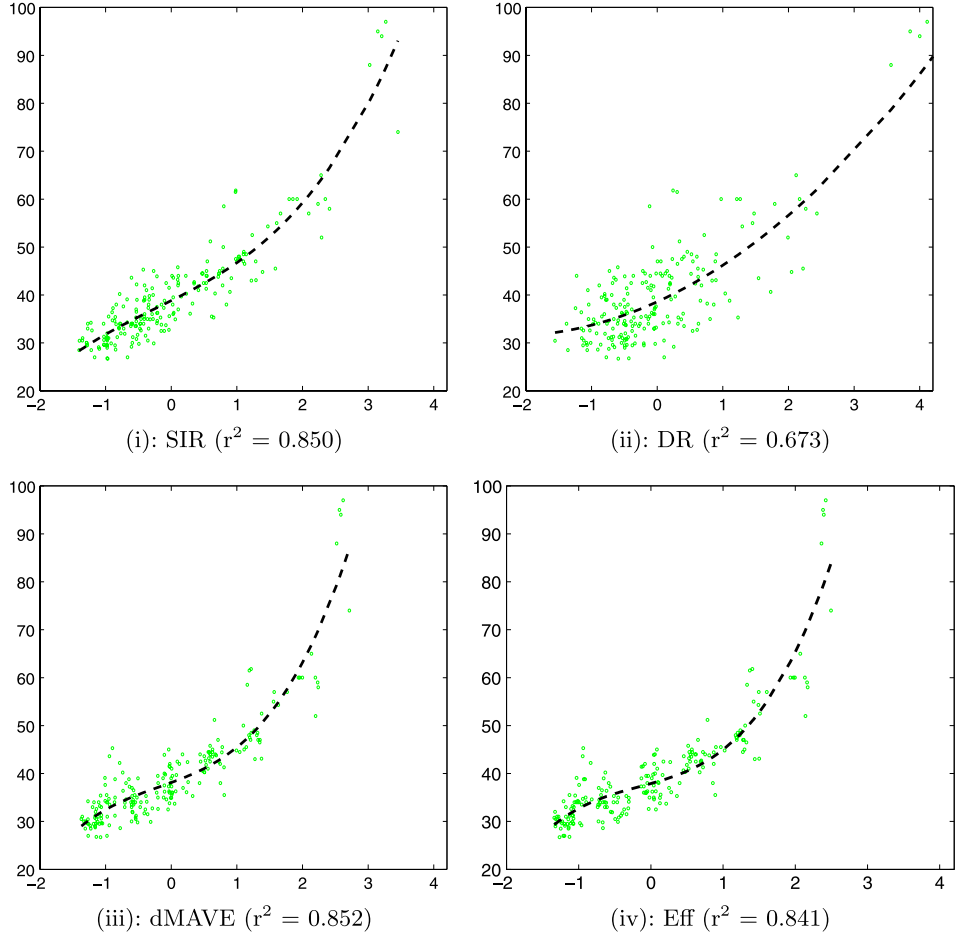}

\caption{The scatter plot of $Y$ versus $\hb^{\mathrm{T}}\x$, with
$\hb$ obtained
from SIR, DR, dMAVE and Eff, respectively.
The fitted cubic regression curves (--) and the adjusted
$\mathrm{r}^2$ values are shown.} \label{figfig1}
\end{figure}

We calculated the Pearson correlation coefficients
and found the current job level
($X_1$) has the largest correlation with his/her annual salary
($Y$) [$\corr(X_1, Y) =0.614$]. This implies that the
current job level is possibly an important factor and thus we
fix the coefficient of $X_1$ to be 1 in
our subsequent analysis. We applied SIR, DR, dMAVE and Eff
methods to estimate the remaining coefficients. In
Figure~\ref{figfig1} we present the scatter plots of $Y$ versus
a single linear combination $\hb^{\mathrm{T}}\x$, where $\x=
(X_1,\ldots,X_6)^{\mathrm{T}}$ and $\hb$ denote the estimate
obtained from
the four estimation procedures. The scatter plots exhibit
similar monotone patterns in that the annual salary increases with
the value of $\hb^{\mathrm{T}}\x$. Except for DR, the
data cloud of all other three proposals looks very compact. To
quantify this visual difference, we fit a cubic model by regressing\vadjust{\goodbreak}
$Y$ on 1, $(\hb^{\mathrm{T}}\x), (\hb^{\mathrm{T}}\x)^2$ and
$(\hb^{\mathrm{T}}\x)^3$. The
adjusted $\mathrm{r}^2$ values are also reported
in Figure~\ref{figfig1}. The $\mathrm{r}^2$
value of DR is much smaller than that of the
other estimators, which suggests worse
performance of DR. This is not a surprise
because DR requires the most stringent
conditions on the covariate vector~$\x$, which are violated here
because of the categorical covariates. The $\mathrm{r}^2$ values of all other
estimators including Eff are satisfactory, indicating that $\cs$ is
possibly one
dimensional. We would also like to point out that because
the $\mathrm{r}^2$ value factors in the goodness-of-fit of the cubic
model, hence it only provides a reference.

%
\begin{table}
\caption{The estimated coefficients and standard
errors obtained by Eff}\label{tabletab4}
\begin{tabular*}{\textwidth}{@{\extracolsep{\fill}}lcccccc@{}}
\hline
&&\multicolumn{1}{c}{$\bolds{\widehat\beta_2}$}&\multicolumn
{1}{c}{$\bolds{\widehat\beta_3}$}&\multicolumn{1}{c}{$\bolds
{\widehat\beta_4}$}&
\multicolumn{1}{c}{$\bolds{\widehat\beta_5}$}& \multicolumn
{1}{c@{}}{$\bolds{\widehat\beta_6}$}\\
\hline
Eff &coef.& 0.477 & 0.265 & 0.024 & 0.050& 0.146 \\
&std. & 0.021 & 0.031 & 0.030 & 0.037& 0.031 \\
& $p$-value &$<\!10^{-4}$ & $<\!10^{-4}$ & 0.427 & 0.176 &$<\!10^{-4}$
\\
\hline
\end{tabular*}
\end{table}

Table~\ref{tabletab4} contains the estimated coefficients $\widehat
\beta_i$'s,
the standard errors and $p$-values obtained through Eff.\vadjust{\goodbreak}
It can be seen that in addition to the current job
level ($X_1$), working experience at the current bank ($X_2$), age ($X_3$)
and whether or not the job is computer related ($X_6$)
are also important factors on salary. While it is not difficult
to understand the importance of most of these factors, we believe the
age effect
is probably caused by its high correlation with the working
experience [$\corr(X_2,X_3) = 0.676$].

\section{Discussion}\label{secdiscuss}

We have derived both locally efficient and efficient estimators
which exhaust the entire central subspace without imposing
any distributional assumptions. We point out here that
if the linearity condition holds, the efficiency bound does not change. However,
the linearity condition will enable a simplification of the
computation because we can simply plug
$E(\x\mid\bb^{\mathrm{T}}\x)=\bb(\bb^{\mathrm{T}}\bb)^{-1}\bb^{\mathrm{T}}\x$
into the estimation
equation instead of estimating it nonparametrically. However, the
constant variance condition does not seem to contribute to
the efficiency bound or to the computational simplicity. It is
therefore a redundant condition in the efficient estimation of the
central subspace.

In this paper we did not discuss how to determine $d$, the structural
dimension of $\cs$ when
an efficient estimation procedure is used, although we agree that this
is an
important issue in the area of dimension reduction. In the real-data
example, we infer the
structural dimension through the adjusted $\mathrm{r}^2$ values.
This seems a reasonable choice, but the turnout may depend on
how to recover the underlying model structure. How to prescribe a
rigorous data-driven procedure is needed in future works.

Various model extensions have been considered in the dimensional reduction
literature. For example, in partial dimension
reduction problems~\cite{ChiaromonteCookLi2002},
it is assumed that
$F(Y\mid\x)=F(Y\mid\bb^{\mathrm{T}}\x_1,\x_2)$. Here, $\x_1$ is a
covariate sub-vector of $\x$ that the dimension reduction procedure
focuses on, while $\x_2$
is a covariate sub-vector that is known to directly enter the model based
on scientific understanding or convention.
We can see that the
semiparametric analysis and the efficient estimation results derived
here can be adapted to these models, through
changing $\bb^{\mathrm{T}}\x$ to
$(\bb^{\mathrm{T}}\x_1, \x_2)$ in all the corresponding functions and
expectations while everything else remains unchanged.
Another extension is the group-wise dimension
reduction~\cite{LiLiZhu2010}, where the model
$E(Y\mid\x)=\sum_{i=1}^k m_i(Y,\x_i^{\mathrm{T}}\bb_i)$ is considered.
The semiparametric analysis in such models
requires separate investigation,
and it will be interesting to study the efficient
estimation.


%

\begin{supplement}[id=suppA]
\stitle{Supplement to ``Efficient estimation in sufficient dimension
reduction''}
\slink[doi]{10.1214/12-AOS1072SUPP} 
\sdatatype{.pdf}
\sfilename{aos1072\_supp.pdf}
\sdescription{The supplement file aos1072\_supp.\break pdf is available upon request.
It contains derivations of the efficient score for model~(\ref{eqcs})
and an outline of proof for Theorems~\ref{thlocal} and~\ref{theff}.}
\end{supplement}

%


\printaddresses


\begin{thebibliography}{31}

\bibitem{AlbrightWinstonZappe1999}
\begin{bbook}[auto:STB|2013/01/23|16:20:06]
\bauthor{\bsnm{Albright},~\bfnm{S.~C.}\binits{S.~C.}},
  \bauthor{\bsnm{Winston},~\bfnm{W.~L.}\binits{W.~L.}} \AND
  \bauthor{\bsnm{Zappe},~\bfnm{C.~J.}\binits{C.~J.}}
(\byear{1999}).
\btitle{Data Analysis and Decision Making with Microsoft Excel}.
\bpublisher{Duxbury}, \blocation{Pacific Grove, CA}.
\bptok{imsref}%
\end{bbook}
\endbibitem

\bibitem{BickelKlaassenRitovWellner1993}
\begin{bbook}[mr]
\bauthor{\bsnm{Bickel},~\bfnm{Peter~J.}\binits{P.~J.}},
  \bauthor{\bsnm{Klaassen},~\bfnm{Chris A.~J.}\binits{C.~A.~J.}},
  \bauthor{\bsnm{Ritov},~\bfnm{Ya'acov}\binits{Y.}} \AND
  \bauthor{\bsnm{Wellner},~\bfnm{Jon~A.}\binits{J.~A.}}
(\byear{1993}).
\btitle{Efficient and Adaptive Estimation for Semiparametric Models}.
\bpublisher{Johns Hopkins Univ. Press}, \blocation{Baltimore, MD}.
\bid{mr={1245941}}
\bptok{imsref}%
\end{bbook}
\endbibitem

\bibitem{ChiaromonteCookLi2002}
\begin{barticle}[mr]
\bauthor{\bsnm{Chiaromonte},~\bfnm{Francesca}\binits{F.}},
  \bauthor{\bsnm{Cook},~\bfnm{R.~Dennis}\binits{R.~D.}} \AND
  \bauthor{\bsnm{Li},~\bfnm{Bing}\binits{B.}}
(\byear{2002}).
\btitle{Sufficient dimension reduction in regressions with categorical
  predictors}.
\bjournal{Ann. Statist.}
\bvolume{30}
\bpages{475--497}.
\bid{doi={10.1214/aos/1021379862}, issn={0090-5364}, mr={1902896}}
\bptok{imsref}%
\end{barticle}
\endbibitem

\bibitem{Cook1994}
\begin{barticle}[mr]
\bauthor{\bsnm{Cook},~\bfnm{R.~Dennis}\binits{R.~D.}}
(\byear{1994}).
\btitle{On the interpretation of regression plots}.
\bjournal{J. Amer. Statist. Assoc.}
\bvolume{89}
\bpages{177--189}.
\bid{issn={0162-1459}, mr={1266295}}
\bptok{imsref}%
\end{barticle}
\endbibitem

\bibitem{Cook1998}
\begin{bbook}[mr]
\bauthor{\bsnm{Cook},~\bfnm{R.~Dennis}\binits{R.~D.}}
(\byear{1998}).
\btitle{Regression Graphics}.
\bpublisher{Wiley}, \blocation{New York}.
\bid{doi={10.1002/9780470316931}, mr={1645673}}
\bptok{imsref}%
\end{bbook}
\endbibitem

\bibitem{CookWeisberg1991}
\begin{barticle}[mr]
\bauthor{\bsnm{Cook},~\bfnm{R.~D.}\binits{R.~D.}} \AND
  \bauthor{\bsnm{Weisberg},~\bfnm{S.}\binits{S.}}
(\byear{1991}).
\btitle{Comment on ``Sliced inverse regression for dimension reduction,'' by K.-C. Li}.
\bjournal{J. Amer. Statist. Assoc.}
\bvolume{86}
\bpages{328--332}.
\bptok{imsref}%
\end{barticle}
\endbibitem

\bibitem{DongLi2010}
\begin{barticle}[mr]
\bauthor{\bsnm{Dong},~\bfnm{Yuexiao}\binits{Y.}} \AND
  \bauthor{\bsnm{Li},~\bfnm{Bing}\binits{B.}}
(\byear{2010}).
\btitle{Dimension reduction for non-elliptically distributed predictors:
  Second-order methods}.
\bjournal{Biometrika}
\bvolume{97}
\bpages{279--294}.
\bid{doi={10.1093/biomet/asq016}, issn={0006-3444}, mr={2650738}}
\bptok{imsref}%
\end{barticle}
\endbibitem

\bibitem{FanYaoTong1996}
\begin{barticle}[mr]
\bauthor{\bsnm{Fan},~\bfnm{Jianqing}\binits{J.}},
  \bauthor{\bsnm{Yao},~\bfnm{Qiwei}\binits{Q.}} \AND
  \bauthor{\bsnm{Tong},~\bfnm{Howell}\binits{H.}}
(\byear{1996}).
\btitle{Estimation of conditional densities and sensitivity measures in
  nonlinear dynamical systems}.
\bjournal{Biometrika}
\bvolume{83}
\bpages{189--206}.
\bid{doi={10.1093/biomet/83.1.189}, issn={0006-3444}, mr={1399164}}
\bptok{imsref}%
\end{barticle}
\endbibitem

\bibitem{LiDong2009}
\begin{barticle}[mr]
\bauthor{\bsnm{Li},~\bfnm{Bing}\binits{B.}} \AND
  \bauthor{\bsnm{Dong},~\bfnm{Yuexiao}\binits{Y.}}
(\byear{2009}).
\btitle{Dimension reduction for nonelliptically distributed predictors}.
\bjournal{Ann. Statist.}
\bvolume{37}
\bpages{1272--1298}.
\bid{doi={10.1214/08-AOS598}, issn={0090-5364}, mr={2509074}}
\bptok{imsref}%
\end{barticle}
\endbibitem

\bibitem{LiWang2007}
\begin{barticle}[mr]
\bauthor{\bsnm{Li},~\bfnm{Bing}\binits{B.}} \AND
  \bauthor{\bsnm{Wang},~\bfnm{Shaoli}\binits{S.}}
(\byear{2007}).
\btitle{On directional regression for dimension reduction}.
\bjournal{J.~Amer. Statist. Assoc.}
\bvolume{102}
\bpages{997--1008}.
\bid{doi={10.1198/016214507000000536}, issn={0162-1459}, mr={2354409}}
\bptok{imsref}%
\end{barticle}
\endbibitem

\bibitem{LiZhaChiaromonte2005}
\begin{barticle}[mr]
\bauthor{\bsnm{Li},~\bfnm{Bing}\binits{B.}},
  \bauthor{\bsnm{Zha},~\bfnm{Hongyuan}\binits{H.}} \AND
  \bauthor{\bsnm{Chiaromonte},~\bfnm{Francesca}\binits{F.}}
(\byear{2005}).
\btitle{Contour regression: A general approach to dimension reduction}.
\bjournal{Ann. Statist.}
\bvolume{33}
\bpages{1580--1616}.
\bid{doi={10.1214/009053605000000192}, issn={0090-5364}, mr={2166556}}
\bptok{imsref}%
\end{barticle}
\endbibitem

\bibitem{Li1991}
\begin{barticle}[mr]
\bauthor{\bsnm{Li},~\bfnm{Ker-Chau}\binits{K.-C.}}
(\byear{1991}).
\btitle{Sliced inverse regression for dimension reduction (with discussion)}.
\bjournal{J. Amer. Statist. Assoc.}
\bvolume{86}
\bpages{316--342}.
\bid{issn={0162-1459}, mr={1137117}}
\bptok{imsref}%
\end{barticle}
\endbibitem

\bibitem{LiDuan1991}
\begin{barticle}[mr]
\bauthor{\bsnm{Li},~\bfnm{Ker-Chau}\binits{K.-C.}} \AND
  \bauthor{\bsnm{Duan},~\bfnm{Naihua}\binits{N.}}
(\byear{1989}).
\btitle{Regression analysis under link violation}.
\bjournal{Ann. Statist.}
\bvolume{17}
\bpages{1009--1052}.
\bid{doi={10.1214/aos/1176347254}, issn={0090-5364}, mr={1015136}}
\bptnote{check year}%
\bptok{imsref}%
\end{barticle}
\endbibitem

\bibitem{LiLiZhu2010}
\begin{barticle}[mr]
\bauthor{\bsnm{Li},~\bfnm{Lexin}\binits{L.}},
  \bauthor{\bsnm{Li},~\bfnm{Bing}\binits{B.}} \AND
  \bauthor{\bsnm{Zhu},~\bfnm{Li-Xing}\binits{L.-X.}}
(\byear{2010}).
\btitle{Groupwise dimension reduction}.
\bjournal{J. Amer. Statist. Assoc.}
\bvolume{105}
\bpages{1188--1201}.
\bid{doi={10.1198/jasa.2010.tm09643}, issn={0162-1459}, mr={2752614}}
\bptok{imsref}%
\end{barticle}
\endbibitem

\bibitem{MaCarroll2006}
\begin{barticle}[mr]
\bauthor{\bsnm{Ma},~\bfnm{Yanyuan}\binits{Y.}} \AND
  \bauthor{\bsnm{Carroll},~\bfnm{Raymond~J.}\binits{R.~J.}}
(\byear{2006}).
\btitle{Locally efficient estimators for semiparametric models with measurement
  error}.
\bjournal{J. Amer. Statist. Assoc.}
\bvolume{101}
\bpages{1465--1474}.
\bid{doi={10.1198/016214506000000519}, issn={0162-1459}, mr={2279472}}
\bptok{imsref}%
\end{barticle}
\endbibitem

\bibitem{MaChiouWang2006}
\begin{barticle}[mr]
\bauthor{\bsnm{Ma},~\bfnm{Yanyuan}\binits{Y.}},
  \bauthor{\bsnm{Chiou},~\bfnm{Jeng-Min}\binits{J.-M.}} \AND
  \bauthor{\bsnm{Wang},~\bfnm{Naisyin}\binits{N.}}
(\byear{2006}).
\btitle{Efficient semiparametric estimator for heteroscedastic partially linear
  models}.
\bjournal{Biometrika}
\bvolume{93}
\bpages{75--84}.
\bid{doi={10.1093/biomet/93.1.75}, issn={0006-3444}, mr={2277741}}
\bptok{imsref}%
\end{barticle}
\endbibitem

\bibitem{MaGenton2010}
\begin{barticle}[mr]
\bauthor{\bsnm{Ma},~\bfnm{Yanyuan}\binits{Y.}} \AND
  \bauthor{\bsnm{Genton},~\bfnm{Marc~G.}\binits{M.~G.}}
(\byear{2010}).
\btitle{Explicit estimating equations for semiparametric generalized linear
  latent variable models}.
\bjournal{J. R. Stat. Soc. Ser. B Stat. Methodol.}
\bvolume{72}
\bpages{475--495}.
\bid{doi={10.1111/j.1467-9868.2010.00741.x}, issn={1369-7412}, mr={2758524}}
\bptok{imsref}%
\end{barticle}
\endbibitem

\bibitem{MaGentonTsiatis2005}
\begin{barticle}[mr]
\bauthor{\bsnm{Ma},~\bfnm{Yanyuan}\binits{Y.}},
  \bauthor{\bsnm{Genton},~\bfnm{Marc~G.}\binits{M.~G.}} \AND
  \bauthor{\bsnm{Tsiatis},~\bfnm{Anastasios~A.}\binits{A.~A.}}
(\byear{2005}).
\btitle{Locally efficient semiparametric estimators for generalized
  skew-elliptical distributions}.
\bjournal{J. Amer. Statist. Assoc.}
\bvolume{100}
\bpages{980--989}.
\bid{doi={10.1198/016214505000000079}, issn={0162-1459}, mr={2201024}}
\bptok{imsref}%
\end{barticle}
\endbibitem

\bibitem{MaHart2007}
\begin{barticle}[mr]
\bauthor{\bsnm{Ma},~\bfnm{Yanyuan}\binits{Y.}} \AND
  \bauthor{\bsnm{Hart},~\bfnm{Jeffrey~D.}\binits{J.~D.}}
(\byear{2007}).
\btitle{Constrained local likelihood estimators for semiparametric skew-normal
  distributions}.
\bjournal{Biometrika}
\bvolume{94}
\bpages{119--134}.
\bid{doi={10.1093/biomet/asm020}, issn={0006-3444}, mr={2307902}}
\bptok{imsref}%
\end{barticle}
\endbibitem

\bibitem{MaZhu2012}
\begin{barticle}[mr]
\bauthor{\bsnm{Ma},~\bfnm{Yanyuan}\binits{Y.}} \AND
  \bauthor{\bsnm{Zhu},~\bfnm{Liping}\binits{L.}}
(\byear{2012}).
\btitle{A semiparametric approach to dimension reduction}.
\bjournal{J.~Amer. Statist. Assoc.}
\bvolume{107}
\bpages{168--179}.
\bid{doi={10.1080/01621459.2011.646925}, issn={0162-1459}, mr={2949349}}
\bptok{imsref}%
\end{barticle}
\endbibitem

\bibitem{supp}
\begin{bmisc}[auto]
\bauthor{\bsnm{Ma},~\bfnm{Yanyuan}\binits{Y.}} \AND
  \bauthor{\bsnm{Zhu},~\bfnm{Liping}\binits{L.}}
(\byear{2013}).
\bhowpublished{Supplement to ``Efficient estimation in sufficient dimension
  reduction.'' DOI:\doiurl{10.1214/12-AOS1072SUPP}.}
\bptok{imsref}%
\end{bmisc}
\endbibitem

\bibitem{Newey1990}
\begin{barticle}[auto:STB|2013/01/23|16:20:06]
\bauthor{\bsnm{Newey},~\bfnm{W.}\binits{W.}}
(\byear{1990}).
\btitle{Semiparametric efficiency bounds}.
\bjournal{J. Appl. Econometrics}
\bvolume{5}
\bpages{99--135}.
\bptok{imsref}%
\end{barticle}
\endbibitem

\bibitem{RobinsRotnitzkyZhao1994}
\begin{barticle}[mr]
\bauthor{\bsnm{Robins},~\bfnm{James~M.}\binits{J.~M.}},
  \bauthor{\bsnm{Rotnitzky},~\bfnm{Andrea}\binits{A.}} \AND
  \bauthor{\bsnm{Zhao},~\bfnm{Lue~Ping}\binits{L.~P.}}
(\byear{1994}).
\btitle{Estimation of regression coefficients when some regressors are not
  always observed}.
\bjournal{J. Amer. Statist. Assoc.}
\bvolume{89}
\bpages{846--866}.
\bid{issn={0162-1459}, mr={1294730}}
\bptok{imsref}%
\end{barticle}
\endbibitem

\bibitem{Tsiatis2006}
\begin{bbook}[mr]
\bauthor{\bsnm{Tsiatis},~\bfnm{Anastasios~A.}\binits{A.~A.}}
(\byear{2006}).
\btitle{Semiparametric Theory and Missing Data}.
\bpublisher{Springer}, \blocation{New York}.
\bid{mr={2233926}}
\bptok{imsref}%
\end{bbook}
\endbibitem

\bibitem{TsiatisMa2004}
\begin{barticle}[mr]
\bauthor{\bsnm{Tsiatis},~\bfnm{Anastasios~A.}\binits{A.~A.}} \AND
  \bauthor{\bsnm{Ma},~\bfnm{Yanyuan}\binits{Y.}}
(\byear{2004}).
\btitle{Locally efficient semiparametric estimators for functional measurement
  error models}.
\bjournal{Biometrika}
\bvolume{91}
\bpages{835--848}.
\bid{doi={10.1093/biomet/91.4.835}, issn={0006-3444}, mr={2126036}}
\bptok{imsref}%
\end{barticle}
\endbibitem

\bibitem{Xia2007}
\begin{barticle}[mr]
\bauthor{\bsnm{Xia},~\bfnm{Yingcun}\binits{Y.}}
(\byear{2007}).
\btitle{A constructive approach to the estimation of dimension reduction
  directions}.
\bjournal{Ann. Statist.}
\bvolume{35}
\bpages{2654--2690}.
\bid{doi={10.1214/009053607000000352}, issn={0090-5364}, mr={2382662}}
\bptok{imsref}%
\end{barticle}
\endbibitem

\bibitem{ZengLin2007a}
\begin{barticle}[auto:STB|2013/01/23|16:20:06]
\bauthor{\bsnm{Zeng},~\bfnm{D.}\binits{D.}} \AND
  \bauthor{\bsnm{Lin},~\bfnm{D.~Y.}\binits{D.~Y.}}
(\byear{2007}).
\btitle{Efficient estimation in the accelerated failure time model}.
\bjournal{J.~Amer. Statist. Assoc.}
\bvolume{102}
\bpages{1387--1396}.
\bptok{imsref}%
\end{barticle}
\endbibitem

\bibitem{ZengLin2007b}
\begin{barticle}[auto:STB|2013/01/23|16:20:06]
\bauthor{\bsnm{Zeng},~\bfnm{D.}\binits{D.}} \AND
  \bauthor{\bsnm{Lin},~\bfnm{D.~Y.}\binits{D.~Y.}}
(\byear{2007}).
\btitle{Maximum likelihood estimation in semiparametric models with censored
  data (with discussion)}.
\bjournal{J. Roy. Statist. Soc. Ser. B}
\bvolume{69}
\bpages{507--564}.
\bptok{imsref}%
\end{barticle}
\endbibitem

\bibitem{ZhuZhuFeng2010}
\begin{barticle}[mr]
\bauthor{\bsnm{Zhu},~\bfnm{Li-Ping}\binits{L.-P.}},
  \bauthor{\bsnm{Zhu},~\bfnm{Li-Xing}\binits{L.-X.}} \AND
  \bauthor{\bsnm{Feng},~\bfnm{Zheng-Hui}\binits{Z.-H.}}
(\byear{2010}).
\btitle{Dimension reduction in regressions through cumulative slicing
  estimation}.
\bjournal{J. Amer. Statist. Assoc.}
\bvolume{105}
\bpages{1455--1466}.
\bid{doi={10.1198/jasa.2010.tm09666}, issn={0162-1459}, mr={2796563}}
\bptok{imsref}%
\end{barticle}
\endbibitem

\bibitem{ZhuZeng2006}
\begin{barticle}[mr]
\bauthor{\bsnm{Zhu},~\bfnm{Yu}\binits{Y.}} \AND
  \bauthor{\bsnm{Zeng},~\bfnm{Peng}\binits{P.}}
(\byear{2006}).
\btitle{Fourier methods for estimating the central subspace and the central
  mean subspace in regression}.
\bjournal{J. Amer. Statist. Assoc.}
\bvolume{101}
\bpages{1638--1651}.
\bid{doi={10.1198/016214506000000140}, issn={0162-1459}, mr={2279485}}
\bptok{imsref}%
\end{barticle}
\endbibitem

\end{thebibliography}
\end{document}